\documentclass[12pt]{amsart}
\usepackage{amscd,amssymb}
\usepackage[arrow,matrix]{xy}
\usepackage[plainpages,backref,urlcolor=blue]{hyperref}

\theoremstyle{plain}
\newtheorem{thm}[subsection]{Theorem}
\newtheorem{lem}[subsection]{Lemma}
\newtheorem{prop}[subsection]{Proposition}
\newtheorem{cor}[subsection]{Corollary}

\theoremstyle{definition}
\newtheorem{rk}[subsection]{Remark}
\newtheorem{definition}[subsection]{Definition}
\newtheorem{ex}[subsection]{Example}
\newtheorem{question}[subsection]{Question}

\numberwithin{equation}{section}
\newcommand{\OO}{{\mathcal O}}

\newcommand{\I}{{\mathcal I}}

\newcommand{\wI}{\widehat{I}}
\newcommand{\wJ}{\widehat{J}}

\newcommand{\NN}{{\mathcal N}}

\newcommand{\Z}{\mathbb{Z}}

\newcommand{\C}{\mathbb{C}}

\newcommand{\PP}{\mathbb{P}}

\newcommand{\N}{\mathbb{N}}

\DeclareMathOperator{\Hom}{Hom}

\DeclareMathOperator{\defect}{def}
\DeclareMathOperator{\codim}{codim}

\DeclareMathOperator{\reg}{reg}

\begin{document}

\title [Syzygies of Jacobian ideals and defects of linear systems]
{Syzygies of Jacobian ideals and defects of linear systems }

\author[Alexandru Dimca]{Alexandru Dimca}
\address{Institut Universitaire de France et Laboratoire J.A. Dieudonn\'e, UMR du CNRS 7351,
                 Universit\'e de Nice Sophia-Antipolis,
                 Parc Valrose,
                 06108 Nice Cedex 02,
                 France}
\email{dimca@unice.fr}



\begin{abstract} 

Our main result  describes the relation between the syzygies involving the first order partial derivatives $f_0,...,f_n$ of a homogeneous polynomial $f\in \C[x_0,...x_n]$ and the defect of the linear systems vanishing on the 
singular locus subscheme $\Sigma_f=V(f_0,...,f_n)$ of the hypersurface $D:f=0$ in the complex projective space $\PP^n$, when $D$ has only isolated singularities. 
\end{abstract}

\maketitle


\section{Introduction} \label{sec:intro}

Let $S=\C[x_0,...,x_n]$ be the graded ring of polynomials in $x_0,,...,x_n$ with complex coefficients and denote by $S_r$ the vector space of homogeneous polynomials in $S$ of degree $r$. 
For any polynomial $f \in S_d$ we define the {\it Jacobian ideal} $J_f \subset S$ as the ideal spanned by the partial derivatives $f_0,...,f_n$ of $f$ with respect to $x_0,...,x_n$. For $n=2$ we use $x,y,z$ instead of
$x_0,x_1,x_2$ and
$f_x,f_y,f_z$  instead of $f_0,f_1,f_2$. 

We define the corresponding graded {\it Milnor} (or {\it Jacobian}) {\it algebra} by
\begin{equation} 
\label{eq1}
M(f)=S/J_f.
\end{equation}
The study of such Milnor algebras is related to the singularities of the corresponding projective hypersurface $D:f=0$, see \cite{CD}, as well as to the mixed Hodge theory of 
the hypersurface $D$ and of its complement $U=\PP^n \setminus D$, see the foundational article by Griffiths \cite{Gr} and also  \cite{DSW}, \cite{DSt2}, \cite{DSt3}, \cite{DS}. 

\bigskip

We define the {\it singular locus scheme} of the hypersurface $D$ to be the subscheme $\Sigma_f$ of $\PP^n$ defined by the ideal $J_f$. If $p$ is an isolated singularity of the hypersurface $D$ with local equation $g=0$, then one has a natural isomorphism
\begin{equation} 
\label{eqT}
\OO_{\Sigma_f,p}=T(g),
\end{equation}
the local Tjurina algebra of the analytic germ $g$, see Lemma \ref{Tj} below. In particular, $\dim \OO_{\Sigma_f,p}=\dim T(g) =\tau(g)$, the Tjurina number of the
isolated singularity $(D,p)$.

\bigskip

On the other hand, one knows that several homogeneous ideals in $S$ may define the same subscheme. The largest one defining $\Sigma_f$
is denoted by $\wJ _f$ and it is the saturated ideal associated to the Jacobian ideal $J_f$, see \cite{H}, p. 125, Exercises II.5.9 and II.5.10. From the definition, it is clear that the homogeneous components $\wJ _{f,k}$ and $J_{f,k}$ of $\wJ _f$ and $J_f$ respectively coincide for $k$ large enough.

In order to get explicit values for such $k$'s, we recall the following notions from \cite{DSt2}.
\begin{definition}
\label{def}
For a degree $d$ hypersurface $D:f=0$ with isolated singularities in $\PP^n$, three integers have been introduced, see \cite{DSt2}.

\noindent (i) the {\it coincidence threshold} $ct(D)$ defined as
$$ct(D)=\max \{q~~:~~\dim M(f)_k=\dim M(f_s)_k \text{ for all } k \leq q\},$$
with $f_s$  a homogeneous polynomial in $S$ of degree $d$ such that $D_s:f_s=0$ is a smooth hypersurface in $\PP^n$.

\noindent (ii) the {\it stability threshold} $st(D)$ defined as
$$st(D)=\min \{q~~:~~\dim M(f)_k=\tau(D) \text{ for all } k \geq q\}$$
where $\tau(D)$ is the total Tjurina number of $D$, i.e. the sum of all the Tjurina numbers of the singularities of $D$.

\noindent (iii) the {\it minimal degree of a nontrivial relation} $mdr(D)$ defined as
$$mdr(D)=\min \{q~~:~~ H^n(K^*(f))_{q+n}\ne 0\}$$
where $K^*(f)$ is the Koszul complex of $f_0,...,f_n$ with the grading defined in  section 3, see also the formula \eqref{relvscoho}.

\end{definition}
It is clear that one has
\begin{equation} 
\label{REL}
ct(D)=mdr(D)+d-2,
\end{equation} 
using the formula \eqref{dif} below, in other words the main invariants are $ct(D)$ and $st(D)$.
By definition, it follows that for any such hypersurface $D$ which is not smooth, we have $d-2 \leq ct(D) \leq T$ and using \cite{CD} we get $st(D) \leq T+1$, where
we set $T=T(n,d)=(n+1)(d-2)$ .

Recall that Hilbert-Poincar\'e series of a graded $S$-module $E=\oplus_kE_k$ of finite type is defined by 
\begin{equation} 
\label{eq21}
HP(E;t)= \sum_{k\geq 0} (\dim E_k)t^k
\end{equation}
and that we have
\begin{equation} 
\label{eq31}
HP(M(f_s);t)=  \frac{(1-t^{d-1})^{n+1}}{(1-t)^{n+1}}.
\end{equation}
In particular,  it follows that $M(f_s)_j=0$ for $j>T$ and $\dim M(f_s)_j=\dim M(f_s)_{T-j}$ for $0 \leq j \leq T$.

\bigskip

In this note we first discuss the relation between the saturation $\wI$ of a homogeneous ideal $I \subset S$ and the condition that a homogeneous polynomial $g \in S_r$ vanishes on the subscheme $V(I)$ of $\PP^n$ defined by $I$.

Our main result is Theorem \ref{linsys} describing the relation between the syzygies involving the partial derivatives $f_0,...,f_n$ and the defect of the linear systems vanishing on the 
 singular locus subscheme $\Sigma_f=V(f_0,...,f_n)=V(J_f)$ of a projective hypersurface $D:f=0$, when $D$ has only isolated singularities. This extends the nodal case treated in Theorem 1.5 in \cite{DSt2} and uses the full power of  the Cayley-Bacharach Theorem as stated in \cite{EGH}, Theorem CB7, i.e. 
the supports of the subschemes $\Gamma'$ and $\Gamma''$ which are residual to each other
might not be disjoint, see Remark \ref{rk3}. For other, more classical relations between syzygies and algebraic geometry we refer the Eisenbud's book \cite{Eis}.

\bigskip

One consequence is the following relation between the above invariants and the saturation $\wJ$ of the Jacobian ideal $J=J_f$.
\begin{cor}
\label{corD}
If the hypersurface $D:f=0$ has only isolated singularities, then $\wJ _k=J_k$ for $k\geq \max(T-ct(D),st(D))$.
\end{cor}
Further relations involving the $a$-invariant $a(M(f))$ and the Castelnuovo-Mumford regularity 
$\reg M(f)$ of the graded algebra (resp. $S$-module ) $M(f)$ are given below.

We would like to thank Laurent Bus\'e who explained to us  alternative proofs of the main result, based on local cohomology and \v Cech complexes and presented here in Remarks \ref{rkBuse}
and \ref{rkBuse2} and who draw our attention on some small errors in a previous version.

\section{Saturation of ideals and defects of linear systems} \label{sec:two}

For any homogeneous ideal $I$ in $S$ we define its saturation $\wI$ as the set of all elements $s \in S$ such that for any $i=0,...,n$ there is a positive integer $m_i$ such that 
$$x_i^{m_i}s \in I,$$
see \cite{H}, p. 125, Exercise II.5.10. It follows that $\wI$ is also a homogeneous ideal in $S$ and moreover $I$ and $\hat I$ define the same subscheme $V(I)=V(\wI) $ of $\PP^n$.

An ideal $I$ is called saturated if $I=\wI$. One has the following alternative definition for a saturated ideal. 

We say that a homogeneous polynomial $g \in S$ vanishes on the scheme $V(I)$ if for any (closed) point $p$ belonging to the support $ |V(I)|$ of our scheme $V(I)$, the germ of regular function induced by $g$ at $p$ (which is defined up to a unit in the local ring $\OO_{\PP^n,p}$) belongs to the ideal sheaf stalk $\I_{V(I),p}$ of the ideal sheaf $\I_{V(I)}$ defining the subscheme $V(I)$.

Then one can easily see that a homogeneous polynomial $h$ in $\wI$ is exactly a homogeneous polynomial vanishing on the subscheme $V(I)=V(\wI) $. Hence, an ideal $I$ is saturated exactly when it contains all the homogeneous polynomials vanishing on the subscheme $V(I)$.

This proves in particular the following, via Theorem 8 (Lasker's Unmixedness Theorem) in \cite{EGH}.

\begin{prop}
\label{prop1}
If the ideal $I$ is a complete intersection, then $I$ is saturated.

\end{prop}
The simpliest, and most well known version of this result, is of course E. Noether's  "AF+BG" Theorem, see for instance \cite{GH} , p. 703.

\begin{rk}
\label{rk1}
When the subscheme $V(I)$ is reduced, then the saturation $\wI$ coincides with the radical ideal $\sqrt I$. For the singular locus $\Sigma_f$, supposed to be 0-dimensional, this happens exactly when $D$ is a nodal hypersurface.

\end{rk}

For any homogeneous ideal $I$ we consider the graded artinian $S$-module
\begin{equation} 
\label{eqSD}
SD(I)=  \frac{\wI}{I},
\end{equation}
called the {\it saturation defect module} of $I$ and the {\it saturation threshold} $sat(I)$ defined as
\begin{equation} 
\label{eqsat}
sat(I)=\min \{q~~:~~\dim I_k= \dim \wI _k \text{ for all } k \geq q\}.
\end{equation}

When $Y=V(I)$ is a 0-dimensional subscheme in $\PP^n$, we introduce the corresponding sequence of defects
\begin{equation} 
\label{eqDEF}
\defect _kY=\dim H^0(Y, \OO_Y)-\dim \frac{S_k}{\wI_k}.
\end{equation}
For the singular locus $Y=\Sigma_f$, if we set $J=J_f$, this becomes
$$\defect _k\Sigma_f=\tau(D)-\dim \frac{S_k}{\wJ_k}.$$
In particular, if $\Sigma_f \ne \emptyset$, then $\defect _0\Sigma_f=\tau(D)-1.$

Moreover, when $D$ is a nodal hypersurface, one clearly has
$h \in \wJ_k$ if and only if $h$ vanishes on the set of nodes $\NN$ of the hypersurface $D$,
i.e. we get exactly the notion used in \cite{DSt2} and \cite{DSt3}.

The module $SD(J_f)$ was already considered  by Pellikaan in \cite{Pe} under the name of Jacobian module (in a local version, and especially when $\wJ_f$ is a radical ideal).

\begin{rk}
\label{rklocalco}
The above objects can be interpreted in terms of local cohomology, see Appendix 1 in \cite{Eis}
for the definition and basic properties of this cohomology. Let $\bf m$ be the maximal ideal $(x_0,...,x_n)$ in $S$. Then one has just from definitions
$$SD(I)=  \frac{\wI}{I}=H^0_{\bf m}(S/I)$$
and also, via Corollary A1.12 in  \cite{Eis},
$$\defect _kY=\dim H^0(Y, \OO_Y)-\dim \frac{S_k}{\wI_k}=\dim H^1_{\bf m}(S/I)_k$$
where the last subscript $k$ indicates the $k$-th homogeneous component.
The $a$-invariant of the graded standard algebra $M(f)$ is given by
\begin{equation} 
\label{ainv}
a(M(f))=\max\{k: H^1_{\bf m}(M(f))_k \ne 0\},
\end{equation}
and the Castelnuovo-Mumford regularity of the graded $S$-module $M(f)$ is given by
\begin{equation} 
\label{reg}
\reg(M(f))=\min\{k: H^0_{\bf m}(M(f))_{>k} =0  \text{ and } H^1_{\bf m}(M(f))_{>k-1} =0   \},
\end{equation}
see \cite{Ch}.

\end{rk}

\section{Defects and syzygies involving the Jacobian ideal} \label{sec3}

Let $f$ be a homogeneous polynomial of degree $d$ in the polynomial ring $S$ and denote by $f_0,...,f_n$ the corresponding partial derivatives.

One can consider the graded $S-$submodule $AR(f) \subset S^{n+1}$ of {\it all relations} involving the $f_j$'s, namely
$$a=(a_0,...,a_n) \in AR(f)_m$$
if and only if  $a_0f_0+a_1f_1+...+a_nf_n=0$.

Inside $AR(f)$ there is the $S-$submodule of {\it Koszul relations} $KR(f)$, called also the 
submodule of {\it trivial relations}, spanned by the relations $t_{ij} \in AR(f)_{d-1}$ for
$0 \leq i <j \leq n$, where $t_{ij}$ has the $i$-th coordinate equal to $f_j$, the $j$-th coordinate equal to $-f_i$ and the other coordinates zero, see relation \eqref{tij} below.

The quotient module $ER(f)=AR(f)/KR(f)$ may be called the module of {\it essential relations}, or non trivial relations, since it tells us which are the relations which we should add to the Koszul relations in order to get all the relations, or syzygies, involving the $f_j$'s.

One has the following description in terms of global polynomial forms on $\C^{n+1}$. If one denotes $\Omega ^j$ the graded $S$-module of such forms of exterior degree $j$, then

\medskip

\noindent(i) $\Omega ^{n+1}$ is a free $S$-module of rank one generated by $\omega=dx_0 \wedge dx_1 \wedge...\wedge dx_n.$

\medskip

\noindent(ii) $\Omega ^{n}$ is a free $S$-module of rank $n+1$ generated by $\omega_j$ for $j=0,...,n$ where $\omega _j$ is given by the same product as $\omega$ but omitting $dx_j$.

\medskip

\noindent(iii) The kernel of the wedge product $df \wedge: \Omega ^{n} \to \Omega ^{n+1}$ can be identified up to a shift in degree to the module $AR(f)$. Indeed, one has to use the formula
$$df \wedge (\sum_{j=0,n}(-1)^ja_j\omega_j)=(\sum_{j=0,n}a_jf_j)\omega.$$

\noindent(iv) $\Omega ^{n-1}$ is a free $S$-module of rank ${n+1 \choose 2}$ generated by $\omega_{ij}$ for $0 \leq i <j \leq n$ where $\omega _{ij}$ is given by the same product as $\omega$ but omitting $dx_i$ and $dx_j$.

\medskip

\noindent(v) The image of the wedge product $df \wedge: \Omega ^{n-1} \to \Omega ^{n}$ can be identified up to a shift in degree to the submodule $KR(f)$. Indeed, one has to use the formula
$$df \wedge \omega_{ij}=f_i\omega_j-f_j\omega_i.$$

In conclusion, it follows that one has
\begin{equation} 
\label{relvscoho}
ER(f)_m=H^n(K^*(f))_{m+n}
\end{equation}
for any $m \in \N$, where $K^*(f)$ is the Koszul complex of $f_0,...,f_n$ with the natural grading $|x_j|=|dx_j|=1$ defined by 
\begin{equation} 
\label{Koszul}
0 \to \Omega^0 \to \Omega^1 \to ...   \to \Omega^{n+1}\to 0
\end{equation}
with all the arrows  given by the wedge product by $df=f_0dx_0+f_1dx_1+...+f_ndx_n$.

Our main result is the following.

\begin{thm}
\label{linsys}
Let $D:f=0$ be a degree $d$ hypersurface in $\PP^n$ having only isolated singularities. If $\Sigma_f$ denotes its singular locus subscheme, then
$$\dim ER(f)_{nd-2n-1-k}=\dim H^{n}(K^*(f))_{nd-n-1-k}= \defect _k\Sigma_f $$
for $0\leq k \leq nd-2n-1$  and $\dim H^{n}(K^*(f))_j=\tau(D)$ for $j\geq n(d-1)$.

In other words, 
$$\dim M(f)_{T-k}=\dim M(f_s)_k+\defect _k\Sigma_f $$
for $0\leq k \leq nd-2n-1$, where $T=T(n,N)=(n+1)(d-2)$. In particular, if  $\Sigma_f \ne \emptyset$, then $\dim M(f)_T=\tau(D)>0$, i.e. $st(D) \leq T$.
\end{thm}
Note that this Theorem determines the dimensions $\dim M(f)_j$ in terms of defects of linear systems for any $j\geq d-1$, i.e. for all $j$ since
the dimensions $\dim M(f)_j=\dim S_j$ for $j<d-1$ are well known. 

\proof

The proof of this result is based on the same idea as the proof of Theorem 1.5 in \cite{DSt2} where the nodal hypersurfaces are treated. 
However, the use of the Cayley-Bacharach Theorem is now more refined, since in the case at hand we deal with non-reduced scheme  $\Sigma_f $.

Let the coordinates on $\PP^n$ be chosen such that the hyperplane $H_0: x_0=0$ is transversal to $D$, then $\Gamma=V(f_1,...,f_n)$ is a 0-dimensional complete intersection contained in the affine space $U_0=\C^n=\PP^n \setminus H_0$. If we use the coordinates $y_1=x_1,....,y_n=x_n$, then the intersection $D_0=D \cap U_0$ is given by the equation $g(y)=0$, where
$g(y)=f(1,y_1,...,y_n)$. The Euler relation for $f$ yields the relation
$$f_0(1,y)+y_1g_1(y)+...+y_ng_n(y)=d \cdot g(y),$$
where $g_j$ denotes the partial derivative of $g$ with respect to $y_j$.
If $0 \in \C^n$ is an isolated singularity of $D_0$, let $\OO_n$ be the  local ring at analytic germs at the origin, $J_g$ the ideal in $\OO_n$ spanned by the partial derivatives. With this notation, one clearly has the following isomorphisms.

\begin{lem}
\label{Tj}
$\OO_{\Gamma,0}=M_g$ where $M_g:=\OO_n/J_g$ is the Milnor algebra of the germ $g$ and $\OO_{\Sigma_f,0}=T_g$ where $T_g:=\OO_n/((g)+J_g)$ is the Tjurina algebra of the germ $g$.
\end{lem}

The support of $\Gamma$ consists of a finite set of points in $\PP^n$, say $p_1,...,p_r$. A part of these points, say $p_j$ for $j=1,...,q$ are  the singularities of $D$, i.e. the points in the support of $\Sigma_f$.

Assume we have a nonzero element in $H^{n}(K^*(f))_{nN-n-1-k}$ for some $0\leq k \leq s,$ with $s=nN-2n-1$. This is the same as having a relation
$$R_m: a_0f_0+a_1f_1+....a_nf_n=0$$
where  $a_j \in S$ are homogeneous of degree
$m=s-k$ and $R_m$ is not a consequence of the relations
\begin{equation} 
\label{tij} 
T_{ij}:f_jf_i-f_if_j=0.
\end{equation} 
This is equivalent to looking at coefficients $a_0$ module the ideal $(f_1,...,f_n)$.

Since $p_j$ is not a singularity for $D$ for $j>q$, it follows that $f_0(p_j)\ne 0$ in this range.
Hence, for $j>q$, the relation $R_m$ implies that the germ of function induced by $a_0$ at $p_j$ (dividing by some homogeneous polynomial $b_j$ of degree $m$ such that $b_j(p_j)\ne 0$) belongs to the ideal defining  $\Gamma$.

At a singular point $p_j$ with $j \leq q$, we get that $a_0f_0 $ belongs to the ideal defining  $\Gamma$. In other words, assuming that $p_j=0$ and using the Euler relation above, we see that the germ induced by $a_0$ at $p_j$ belongs to the annihilator ideal $Ann(g) \subset M_g$ of the class of $g$ in the local Milnor algebra $M_g$.

We apply now the Cayley-Bacharach Theorem as stated in \cite{EGH}, Theorem CB7, where one should replace 'family of curves' by 'family of hypersurfaces' in the last phrase.

Let $\Gamma'$ and $\Gamma''$ be subscheme of $\Gamma$, residual to one another in $\Gamma$, and such that:

(i) the support of $\Gamma'$ is  contained in the support of $\Gamma$,  at a point in the  set $\{p_{q+1},...,p_r\}$ these two schemes coincides, and at a point $p$ in the set  $\{p_{1},...,p_q\}$ the subscheme $\Gamma'$ of $\Gamma$ is defined by the ideal $Ann(g_p)$ in the local Milnor algebra $M_{g_p}$, where $g_p=0$ is a local equation of $D$ at $p$ and we use the identification given by Lemma \ref{Tj} above.

(ii) the support of $\Gamma''$ is the set  $\{p_{1},...,p_q\}$ and the corresponding local ring at a point $p$ in this set is the local Tjurina algebra $T_{g_p}$. In other words, the corresponding ideal is exactly the principal ideal $(g_p)$ of the Milnor algebra $M_{g_p}$. Recall also that these two ideals $(g_p)$ and $Ann(g_p)$ are orthogonal complements to each other via a nondegenerate pairing on the Gorenstein local ring $M_{g_p}$, see \cite{EGH}. More precisely, we have the following general result, perhaps well known to specialists.

\begin{lem}
\label{Ann}
Let $(A,m)$ be a local Artinian Gorenstein ring containing a field $K$. Then, for any ideal $I \subset A$, one has $Ann(I) =I^{\perp}$, where the orthogonal complement is taken with respect to the  nondegenerate pairing $Q:A \times A \to K$.
\end{lem}

\proof
In fact, one has $K=A/m$ and there is a positive integer $s$ such that $m^s$ has length one, i.e. $m^s = K$. For any $K$-linear map $\rho:A \to K$ inducing an isomorphism on $m^s$, the nondegenerate pairing in the statement above may be given as the composition $Q: A \times A \to A \to K$, where the first arrow is the multiplication in $A$ and the second arrow is $\rho$.

From this construction, it is clear that $Ann(I) \subset I^{\perp}$. Conversely, let $a \in I^{\perp}$, such that we have  $Q(ai)=0$ for any $i \in I$. Suppose there is an $i_0 \in I$ such that $ai_0 \ne0$.
Then there is an element $g \in A$ such that $ai_0g \ne 0$ but $ai_0g \in m^s$, see the discussion on page 312 of the paper \cite{EGH} . This is a contradiction, since it implies $Q(a,i_0g) \ne 0$ or we have $i_0g \in I.$

\endproof

The above discussion implies that the dimension of the family of hypersurfaces $a_0$ of degree $m=s-k$ containing $\Gamma'$ (modulo those containing all of $\Gamma$, which are in fact exactly the elements of the ideal $(f_1,...,f_n)$ in view of Proposition \ref{prop1}) is exactly the dimension of $H^{n}(K^*(f))_{nN-n-1-k}$.

On the other hand, for $s$ as above and $0\leq k\leq s$, the Cayley-Bacharach Theorem says that this dimension is equal to the defect $\defect_k(\Sigma_f)$, thus proving the first claim in Theorem \ref{linsys}.

Next we have 
\begin{equation} 
\label{dif} 
\dim H^{n}(K^*(f))_j
=\dim M(f)_{j+d-n-1}-\dim M(f_s)_{j+d-n-1},
\end{equation} 
see \cite{DSt2}.
Moreover, $j \geq n(d-1)$ is equivalent to $j+d-n-1>(n+1)(d-2)$ and hence $\dim M(f)_{j+d-n-1}=\tau(D)$ and $\dim M(f_s)_{j+d-n-1}=0$, thus proving the second claim in Theorem \ref{linsys}.

\endproof

\begin{rk}
\label{rk3}
It follows by the proof above that a point $p$ is in the support of both subschemes $\Gamma'$ and $\Gamma''$ if and only if the singularity $(D,p)$ is not weighted homogeneous.
Indeed, by K. Saito's result \cite{KS}, this is equivalent to the local equation $g_p=0$ satisfying $g \notin J_{g_p}$.
\end{rk}

\begin{rk}
\label{rkBuse}
One can obtain an alternative proof for Theorem \ref{linsys} as follows.
We construct a double complex $K^{*,*}$ in the following way. The $0$-th line is just the Koszul complex considered in \eqref{Koszul}, but shifted in such a way that the differentials become homogeneous of degree $0$ and $K^{0,0}$=S. In terms of the grading ring $S$, one has
$$K^{p,0}=\wedge^pS^{n+1}(p(d-1)).$$
Then  we build the $p$-th column by replacing $K^{p,0}$ by its \v Cech complex as defined for instance in section A1B of \cite{Eis} or in  \cite{Buse}, p. 18, whose notation we use below. In other words, we set
$$K^{p,q}=C^q({\bf x},\wedge^pS^{n+1}(p(d-1))),$$
where ${\bf x}=(x_0,...,x_n)$.
To get our result one has to consider the associated double complex 
$T^r=\oplus_{p+q=r} K^{p,q}$ and to compute its cohomology in two ways, using the two usual spectral sequences, exactly as in the proof of Lemma 3.13 in  \cite{Buse}.

If we compute first the cohomology along the columns using Theorem A1.3 and Theorem A2.50 in \cite{Eis},
the only non-trivial groups are on the $(n+1)$-st line, and this can be identified to the dual of the Koszul complex. More precisely, we have isomorphism of vector spaces
$$H^{n+1}_{\bf m}(K^{p,0})_s= \wedge ^{n+1-p}S_{(n+1-p)(d-1)-s}=\Hom(K^{n+1-p,0},\C)_s.$$
Next we compute first the cohomology along the lines, and we get nonzero terms only  on the last two columns, which correspond to the \v Cech complex for $H^n(K^*(f))$ (resp. $H^{n+1}(K^*(f))$)
on the $n$-th column (resp. on the $(n+1)$-st column). Both of these $S$-modules have a support of dimension 1, hence when we take now the cohomology along the columns, we get
$$H^{j}_{\bf m}(H^n(K^*(f)))=H^{j}_{\bf m}(H^{n+1}(K^*(f)))=0$$
for $j>1$. On the other hand one clearly has $H^{0}_{\bf m}(H^n(K^*(f)))=0$, see for instance 
\cite{CD}, Corollary 11, and hence $E_2=E_{\infty}$ for this spectral sequence as well.

Putting everything together we get that
$$\dim H^{n+2}(T^*)_s=\dim ER(f)_{nd-n-s}=\dim H^1_{\bf m}(M(f))_{-n-1+s}$$
which is exactly the claim of Theorem \ref{linsys} .

Note that exactly the same proof works for any collection of $n+1$ homogeneous polynomials of the same degree $(d-1)$ when they define a zero-dimensional subscheme of $\PP^n$.
The case of homogeneous polynomials of different degrees can be handled in a similar way, but more care is needed with the homogeneity shifts to assure that the corresponding Koszul complex has degree 0 differentials. See also \cite{vS}.

\end{rk}

\begin{rk}
\label{rkBuse2}
A more rapid proof, essentialy equivalent to the above, can be obtained as follows. First we use the local duality, namely if 
$\omega_{M(f)}$ is the canonical module of $M(f)$, then the dual
$H^1_{\bf m}(M(f))^{\vee}$ is graded isomorphic to $\omega_{M(f)}$, see Theorem 3.6.19 page 142 in \cite{BH} or Fact 5 in \cite{Ch}, where the graded version is clearly stated. Then recall that the first nonzero cohomology group in the Koszul complex is nothing else but the shifted canonical module, namely with the grading for $K^*(f)$ considered in the Remark above and in \cite{Ch} one has
$$(\omega_{M(f)})_j=H^n(K^*(f)) _{j-n-1}=ER(f)_{nd-2n-1+j},$$
see Theorem 1.6.16 page 50 in \cite{BH} and Lemma 22 in \cite{Ch}.

\end{rk}

\begin{ex}
\label{ex3A2}

\textbf {Quartic curves  $(d=4)$} 
Any quartic curve with 3 cusps is projectively isomorphic
to the  curve:
$$ C : f=x^2y^2+y^2z^2+z^2x^2-2xyz(x+y+z)= 0$$ 
In general, for any   cuspidal curve $C$, one can imagine the singular locus subscheme as consisting of a family of points $p$ (located at the cusps of $C$) and a nonzero cotangent vector $u_p$ at every such point $p$ (given by the corresponding tangent cone).
Then a homogeneous polynomial $g$ vanishes on $\Sigma_f$ if and only if one has $g(p)=0$ and $dg(p)=\lambda_pu_p$ for some constants $\lambda_p\in \C$.
A direct computation shows that for our quartic curve above we have
$$HP(M(f);t)=1+3t+6t^{2}+7t^{3}+6(t^{4}+\ldots \\$$
and
$$HP(M(f_s);t)=1+3t+6t^{2}+7t^{3}+6t^{4}+3t^5+t^6.$$
The $3$ cusps are located at the points $a=(0:0:1)$, $b=(0:1:0)$ and $c=(1:0:0)$ and the nonzero cotangent vectors are $u_a=dx-dy$ and so on.
Using Theorem \ref{linsys} we get the following.

\medskip

\noindent (i) $\defect _0\Sigma_f =6-1=5$, $\defect _1\Sigma_f =6-3=3$, $\defect _k\Sigma_f =0$ for $k \geq 2$. Using the definition of $\defect _k\Sigma_f $, this yields $\wJ _k=0$ for $k=0,1,2$, 
(which is clear using our geometric description) and $\dim \wJ_ m={m+2 \choose 2}-6$  for $m \geq 3$ where $J=J_f$.

\medskip

\noindent (ii) On the other hand, we obviously have $J _k=0$ for $k=0,1,2$ and 
$\dim J_m={m+2 \choose 2}-\dim M(f)_m$  for $m \geq 3$.

\medskip

It follows that $sat(J_f)=4=st(D)$ and $SD(J_f)=\C$ placed in degree $3$.

\end{ex}

\section{Some consequences} \label{sec4}

Using the above notations, we have the following.

\begin{prop}
\label{prop2}
Assume the hypersurface $D:f=0$ in $\PP^n$ has only isolated singularities and set $J=J_f$. Then the sequence of dimensions
$\dim \frac{S_k}{\wJ _k}$ is an increasing sequence bounded by the total Tjurina number  $\tau(D)$ of $D$ given by
$$\tau(D)= \sum_{p \in |\Sigma_f|}\tau(g_p).$$
Moreover, $\dim \frac{S_k}{\wJ _k}=\tau(D)$ if and only if $k \geq T-ct(D)$.

\end{prop}

\proof The first claim follows from Theorem \ref{linsys} and Corollary 11 in \cite{CD} which show that the sequence of
defects $\defect_k\Sigma_f$ is decreasing. The second claim follows from the equality \eqref{dif} and the definition of $ct(D)$.

\endproof

\begin{cor}
\label{corB}
$$sat(J_f)\leq max(T-ct(D),st(D)).$$
\end{cor}
The example $f=xyz$ where $sat(J_f)=0$, $T-ct(D)=st(D)=1$ shows that this inequality may be strict.
In fact, based on empirical evidence as seen in the following Example, one may conjecture that $T-ct(D)\leq st(D)$. Note that $J_f=\wJ_f$ implies $T-ct(D)= st(D)$.

\begin{ex}
\label{exnodal}

\noindent (i) Let $D:x^py^q+z^d=0$ where $p>0,$ $q>0$ and $p+q=d$. Then using the partial derivatives $f_x$ and $f_y$ we see that $mdr(D)=1$, and hence $ct(D)=d-1<\frac{T}{2}$.
A direct computation using Example 14. (i) in \cite{CD} yields $st(D)= 2d-3$.

\medskip

\noindent (ii) Let $D$ be a degree $d$ nodal hypersurface in $\PP^n$. Then  $ct(D) \geq \frac{T}{2}$,
see \cite{DSt3}, Corollary 2.2. On the other hand, it is clear that one has in general $st(D) \geq ct(D)$, except possibly the case
when $T$ is odd and $st(D)=ct(D)-1=\frac{T-1}{2}.$ However, note that in this very special case one has $T-ct(D)=st(D)$.

\medskip

Hence  in case (i) as well as for all hypersurfaces $D$ such that $ct(D) \geq \frac{T}{2}$ (as  in case (ii) above),
 we get $T -ct(D) \leq st(D)$ and hence $sat(J_f)\leq st(D)$.
\end{ex}

\begin{prop}
\label{prop3}
Assume the hypersurface $D:f=0$ in $\PP^n$ has only isolated singularities and assume that $st(D) \geq n(d-2)+1=T-(d-3).$ Then $sat(D)=st(D)$.

\end{prop}

\proof Corollary 8 in \cite{CD} shows that $\dim M(f)_{q-1} \geq \dim M(f)_{q}$ for $q\geq n(d-2)+1$. It follows that for $q=st(D)$ one has
$$\codim J_{q-1}=\dim M(f)_{q-1}>\tau(D) \geq \codim\wJ_{q-1}.$$

\endproof

In the following example we list the few general situations where the explicit value of $sat(D)=st(D)$ is known.
\begin{ex}
\label{st}

\noindent (i) Let $D$ be a degree $d$ nodal curve in $\PP^2$, which is not a line arrangement. Then
one has $st(D)\geq 2d-3$, see formula (1.6) and Corollary 1.4 in \cite{DSt2}.
In particular, if $D$ has just one node and $d>2$, then we have $sat(D)=T=3d-6$,
see Example 4.3 (i) in \cite{DSt2}.

\medskip

\noindent (ii) Let $D$ be a degree $d$ Chebyshev hypersurface in $\PP^n$. Then
one has $st(D)=T-(d-3)$, see Corollary 3.2 in \cite{DSt3}.

\end{ex}

\begin{cor}
\label{corC}
Assume that $ct(D) \geq \frac{T}{2}$. Then $\tau(D) \leq \dim M(f_s)_{T-ct(D)}$.
\end{cor}
This results shows that for large  $ct(D)$, i.e. $ct(D)$ close to $T$, the Tjurina number (and in particular the number of singularities) has to be small. For instance, $ct(D)=T$ if and only if  $D$ has only one singularity, and this is of type $A_1$, i.e. a node.

One has also the following result, using the definitions given in \eqref{ainv} and \eqref{reg} 
and Theorem \ref{linsys}.

\begin{cor}
\label{corD}
Assume the hypersurface $D:f=0$ in $\PP^n$ has only isolated singularities and $d=deg(f)$. Then 
$$a(M(f))=nd-2n-1-mdr(D)=T-ct(D)-1$$
and
$$\reg(M(f))=\max(T-ct(D),sat(J_f)-1).$$
\end{cor}

\begin{rk}
\label{sym}
It is shown in \cite{vS} in the local case and in \cite{DS} in the graded case that the torsion module
$SD(J_f)=H^0_{\bf m}(M(f))$ is a Gorenstein module, and hence in particular has interesting symmetry properties. In the graded case this can be stated as
$$\dim SD(J_f)_k=\dim SD(J_f)_{T-k}$$
for all $k \in \Z$. Moreover, we {\it conjecture} that the sequence of dimensions $\dim SD(J_f)_k$ is  {\it unimodal}, i.e. one has $\dim SD(J_f)_k \leq \dim SD(J_f)_{k+1}$ for all $0 \leq k <T/2$.
As an example, when $f=x(x^3+y^3+z^3)$, i.e. $n=2$ and $d=4$, the corresponding sequence is $0,1,3,4,3,1,0$.

\end{rk}
\section{The case $\Sigma_f$ is a complete intersection} \label{sec5}

In this section we show how Theorem \ref{linsys} can be used to obtain a new proof of the following result obtained in \cite{CD}, Proposition 13. Assume as above that the hypersurface $D:f=0$ in $\PP^n$ has only isolated singularities and $d=deg(f)$.

Assume moreover that $\Sigma_f$ is a complete intersection, i.e. there are homogeneous polynomials $g_1,...,g_n$ in $S$ of degrees $a_1,...,a_n$, such that the ideal $I$ in $S$ spanned by the $g_i$'s satisfies the conditions
$$J_f \subset I \text{ and } (J_f)_s=I_s \text{ for all } s>>0.$$
It it clear that this condition can be restated as $J_f \subset I \subset \wJ_f.$ But this implies
$ I \subset \wJ_f \subset \wI$ and using Proposition \ref{prop1} we get
$$I = \wJ_f = \wI,$$
hence the ideal $I$ is precisely the saturation $ \wJ_f $.

With this notation we have the following result.
\begin{prop}
\label{prop4}
$$HP(M(f))(t)=HP(M(f_s))+t^{(n+1)(d-1)-\sum a_i} \frac{(1-t^{a_1}) \cdots (1-t^{a_n})}{(1-t)^{n+1}}.$$

\end{prop}

\proof It follows from Theorem \ref{linsys} that one has
$$\dim M(f)_{k}=\dim M(f_s)_k+\tau(D) -m''_{T-k}, $$
where $m''_j=\dim (S/I)_j$. The Hilbert-Poincar\'e series of $S/I$ is just
$$\frac{(1-t^{a_1}) \cdots (1-t^{a_n})}{(1-t)^{n+1}}.$$
This implies $m''_j+m''_{q-k-1}=\tau(D)$ for all $j \in \Z$, where $q=\sum a_i -n$.
We get
$$\dim M(f)_{k}=\dim M(f_s)_k+m''_{\sum a_i-(n+1)(d-1)+k}, $$
which is equivalent to our claim.

\endproof
This proof yields also the following.
\begin{cor}
\label{corE}
With the notation and assumptions of this section, we have $\tau(D)=a_1\cdots a_n$ and $ct(D)=T-\sum a_i +n$. In particular, when $n=2$, the couple $(a_1,a_2)$ is determined, when it exists, by the couple $(\tau(D),ct(D))$.
\end{cor}

\end{document}